\documentclass[11pt]{article}
\usepackage[utf8]{inputenc}
\usepackage[T1]{fontenc}
\usepackage[a4paper,bindingoffset=0cm,head=15pt,body={17cm,25cm},top=2.5cm]{geometry}
\usepackage{amsmath,amssymb}
\usepackage{graphicx}
\usepackage{color}
\usepackage{cite}
\usepackage{ifthen}
\usepackage{titlesec}
\usepackage{fancyhdr}
\usepackage{booktabs}
\usepackage[english]{babel}
\usepackage{hyperref}
\usepackage{titlecaps}
\usepackage{xparse}
\usepackage{dblfloatfix}

\usepackage{placeins}
\usepackage{tikz}
\usetikzlibrary{arrows}
\usetikzlibrary[patterns]
\usetikzlibrary{shapes.geometric,decorations.markings}
\usetikzlibrary{arrows.meta,decorations.pathreplacing}

\definecolor{fftttt}{rgb}{1,0.2,0.2}
\definecolor{qqqqcc}{rgb}{0,0,0.8}
\definecolor{wwwwww}{rgb}{0.4,0.4,0.4}

\titleformat{\section}{\normalfont\bfseries}{\thesection}{2ex}{}
\titlespacing*{\section}{0pt}{2.5ex}{1ex}
\titleformat{\subsection}{\normalfont\bfseries}{\thesubsection}{2ex}{}
\titlespacing*{\subsection}{0pt}{2.5ex}{1ex}
\newcounter{authorcounter}
\newcounter{institutecounter}

\newcommand{\wavesauthorlist}{}
\newcommand{\wavesaddresslist}{}
\newcommand{\wavesemail}{}
\newcommand{\wavesfootnotes}{}
\newcommand{\wavesauthorpre}{}

\def\theNumberTest#1{%
  \if\relax\detokenize\expandafter{\romannumeral-0#1}\relax
    true%
  \else
    false%
  \fi
}

\NewDocumentCommand{\wavesspeaker}{ O{} O{} m m}{%
    \ifthenelse{\value{authorcounter} > 1}{%
      \renewcommand{\wavesauthorpre}{, }%
    }{%
      \renewcommand{\wavesauthorpre}{}%
    }%
    \ifx\relax#1\relax
      \renewcommand{\wavesfootnotes}{}
    \else   
      \renewcommand{\wavesemail}{$^\ast$Email: #1}%
      \renewcommand{\wavesfootnotes}{, \ast}
    \fi
    \ifthenelse{\equal{\theNumberTest{#4}}{true}}{%
      \edef\wavesauthorlist{\wavesauthorlist%
        \wavesauthorpre{}\underline{#3}$^{#2%
        }$%
      }%
    }{%
      \edef\wavesauthorlist{\wavesauthorlist%
        \wavesauthorpre\underline{#3}$^{\arabic{authorcounter}%
        }$%
      }
      \edef\wavesaddresslist{\wavesaddresslist%
        \par%
        $^{\arabic{authorcounter}}$#4%
      }%
      \stepcounter{authorcounter}%
      \stepcounter{institutecounter}
    }%

    \ifx\relax#2\relax
          \edef\wavesauthorlist{\wavesauthorlist%
        \wavesauthorpre$^{\wavesfootnotes%
        }$%
      }
    \else
    \ifthenelse{\equal{\theNumberTest{#2}}{true}}{%
      \edef\wavesauthorlist{\wavesauthorlist%
        \wavesauthorpre{}$^{,#2\wavesfootnotes%
        }$%
      }%
    }{%
      \edef\wavesauthorlist{\wavesauthorlist%
        \wavesauthorpre$^{,\arabic{institutecounter}\wavesfootnotes%
        }$%
      }
      \edef\wavesaddresslist{\wavesaddresslist%
        \par%
        $^{\arabic{institutecounter}}$#2%
      }%
      \stepcounter{institutecounter}%
    }%
    \fi
  \ignorespaces
}

\NewDocumentCommand{\wavesauthor}{ O{} O{} m m}{%
    \ifthenelse{\value{authorcounter} > 1}{%
      \renewcommand{\wavesauthorpre}{, }%
    }{%
      \renewcommand{\wavesauthorpre}{}%
    }%
    \ifx\relax#1\relax
      \renewcommand{\wavesfootnotes}{}
    \else   
      \renewcommand{\wavesemail}{$^\ast$Email: #1}%
      \renewcommand{\wavesfootnotes}{, \ast}
    \fi%
    \ifthenelse{\equal{\theNumberTest{#4}}{true}}{%
      \edef\wavesauthorlist{\wavesauthorlist%
        \wavesauthorpre{}#3$^{#4%
        }$%
      }%
    }{%
      \edef\wavesauthorlist{\wavesauthorlist%
         \wavesauthorpre{}#3$^{\arabic{institutecounter}%
         }$%
      }
      \edef\wavesaddresslist{\wavesaddresslist%
        \par%
        $^{\arabic{institutecounter}}$#4%
      }%
      \stepcounter{authorcounter}
      \stepcounter{institutecounter}%
    }%

    \ifx\relax#2\relax
          \edef\wavesauthorlist{\wavesauthorlist%
        $^{\wavesfootnotes%
        }$%
      }
    \else
    \ifthenelse{\equal{\theNumberTest{#2}}{true}}{%
      \edef\wavesauthorlist{\wavesauthorlist%
        $^{,#4\wavesfootnotes%
        }$%
      }%
    }{%
      \edef\wavesauthorlist{\wavesauthorlist%
        $^{,\arabic{institutecounter}\wavesfootnotes%
        }$%
      }
      \edef\wavesaddresslist{\wavesaddresslist%
        \par%
        $^{\arabic{institutecounter}}$#2%
      }%
      \stepcounter{institutecounter}%
    }%
    \fi
  \ignorespaces
}

\fancypagestyle{plain}{%
  \fancyhead{}
}

\newenvironment{wavespaper}[3]{%
  \renewcommand{\wavesauthorlist}{}%
  \renewcommand{\wavesemail}{}%
  \setcounter{authorcounter}{1}%
  \setcounter{institutecounter}{1}%
     #2
  \twocolumn[
    \begin{center}
     \bfseries
     #1
     \bigskip

     \wavesauthorlist
     \mdseries
     \smallskip

     \wavesaddresslist
     \smallskip
 
     \wavesemail
    \end{center}%
  ]
}{%
}

\usepackage{mathtools}
\usepackage{subfig}

\begin{document}

\begin{wavespaper}{%

A note on the discrete Unmapped Tent Pitching for the heterogeneous wave equation

}{%
  \wavesspeaker[marcella.bonazzoli@inria.fr]{Marcella Bonazzoli}{Inria, Unit\'e de Math\'ematiques Applique\'es, ENSTA, Institut Polytechnique de Paris, 91120 Palaiseau, France}
  \wavesauthor{Gabriele Ciaramella}{MOX, Dipartimento di Matematica, Politecnico di Milano, 20133, Milano, Italy}
  \wavesauthor[][]{Ilario Mazzieri}{2}
}

\section*{Abstract}
The Unmapped Tent Pitching (UTP) algorithm is a space–time domain decomposition method for the parallel solution of wave-type problems. We have recently extended UTP to heterogeneous settings and compared, at the continuous level, the computational cost of different space–time decompositions. In this note, we report discrete-level observations showing that the optimal decomposition strategy may differ from the one predicted by the continuous analysis.

\smallskip

\noindent\textbf{Keywords:} domain decomposition, space-time solvers, wave equation. 

\section{Introduction}
The UTP algorithm was designed for the parallel solution of the homogeneous wave equation in \cite{procsUTP}. 
It is inspired by the Mapped Tent Pitching (MTP) algorithm \cite{MTP2017}. 
While in MTP each physical tent is mapped onto a space-time Cartesian box in nD (rectangle in 1D) where local problems are solved, UTP avoids the mapping by computing the solution directly on a physical space-time box containing the tent.  
In \cite{UTPhet}, we extended UTP to the heterogeneous case:  
\begin{equation*}
  \partial_{tt} u(x,t) = c(x)^2 \partial_{xx} u(x,t), \, (x,t) \in \Omega\times(0,T),
\end{equation*} 
coupled with homogeneous Dirichlet conditions on $\{0,L\} \times (0,T)$ and non-null initial conditions. 
Here, \(\Omega = (0,L)\), \(T>0\), \(c(x) = c_1 > 0\) for \(x \in (0,L/2]\) and \(c(x) = c_2 > 0\) for \(x \in (L/2,L)\). 
For simplicity, we take \(c_1 > c_2\) and \(T= L/(2c_1)\).

\section{UTP algorithm and discrete setting} 

\begin{figure}[bt]
 \includegraphics[width=0.5\linewidth]{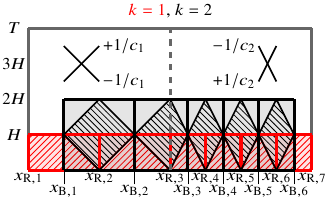} \hspace*{-0.2cm}
 \includegraphics[width=0.5\linewidth]{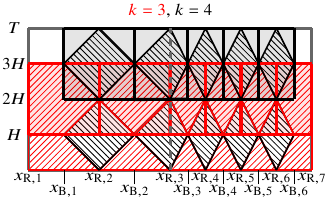}\\
 \caption{First 4 iterations for \(4=m_2>m_1=2\) with \(m_2/m_1 = L_1/L_2 = c_1/c_2\). Red and black boxes denote space–time subdomains; hatched regions indicate exact solution areas.}
 \label{fig:case3}
\end{figure}

UTP is based on a \emph{red-black} decomposition of \(\Omega\times(0,T)\) into space-time subdomains.
On \(\Omega\) we set a grid $\{x_{\mathrm{R},j}\}_{j}$ having $m_1$ and $m_2$ uniformly distributed points in \((0,L/2)\) and \((L/2,L)\), resp. Thus, the \emph{red} subintervals $I_{\mathrm{R},j} \coloneq(x_{\mathrm{R},j},x_{\mathrm{R},j+1})$ have length \(L_1 = L/(2m_1)\) and \(L_2 = L/(2m_2)\) in \((0,L/2)\) and \((L/2,L)\), resp. The \emph{black} subintervals \(I_{\mathrm{B},j} \coloneq (x_{\mathrm{B},j},x_{\mathrm{B},j+1}) \), \(j=1,\dots,m_1+m_2-1\), are defined with \(x_{\mathrm{B},j} \coloneq x_{\mathrm{R},j} + L_1/2\), \(j=1,\dots,m_1\), \(x_{\mathrm{B},j} = x_{\mathrm{R},j} + L_2/2\), \(j=m_1+1,\dots,m_1+m_2\). 
While \cite{UTPhet} also considers strategies with different time heights $H_1$ and $H_2$ in the two regions, here we assume $H_1 = H_2 = H$, with $H = L_1/(2c_1)$ determined by the smaller characteristic slope $\lvert \pm 1/c_1 \rvert$.
More precisely, if \(k\) is the iteration index, we define the time subintervals \(\tau_1 \coloneq (0,H)\) and \(\tau_k \coloneq (t_{k,a}, t_{k,b}) \coloneq \bigl( (k-2)H, \mathrm{min}(kH, T) \bigr) \) for $k\geq2$.
Space-time rectangles \(\mathcal{T}_j^k\) are pitched alternately on the red and black subintervals:  
\(\mathcal{T}_j^k \coloneq I_k \times \tau_k\), with \(I_k \coloneq I_{\mathrm{R},j}\) for \(k\) odd, \(I_k \coloneq I_{\mathrm{B},j}\) for \(k\) even. See, e.g., Fig.~\ref{fig:case3}, where \(c_1=2c_2\). 
Thus, given an approximation $u^{k-1}$ in $\Omega \times (0,T)$ satisfying initial and boundary conditions, the next $u^k$ is computed by solving the subproblems 
\begin{align*}
    &\partial_{tt} u_j^k(x,t) = c(x)^2 \partial_{xx} u_j^k(x,t), \, (x,t) \in \mathcal{T}_j^k, \\
    &u_j^k = u^{k-1}, \, \text{on } \partial I_k \times \tau_k,\\
    &u_j^k = u^{k-1}, \, \partial_t u_j^k = \partial_t u^{k-1}, \, \text{on }  I_k \times \{ t_{k,a} \},
\end{align*}
and then setting \(u^k = u_j^k\) in \(\mathcal{T}_j^k\) for the solved \(j\) subdomains and \(u^k = u_j^{k-1}\) elsewhere. 
The algorithm terminates when the exact solution is computed in the entire \(\Omega\times(0,T)\); see the hatched regions in Fig.~\ref{fig:case3}. 

We discretize $\overline{\Omega}\times[0,T]$ by a uniform time grid of size $\Delta t$, and a piecewise-uniform space grid of size $h_1$ in $[0,L/2]$ and $h_2$ in $[L/2,L]$. The discrete space–time subdomains $\,^h\mathcal{T}_j^k$ are defined on this grid; see Fig.~\ref{fig:discrete_subs}.
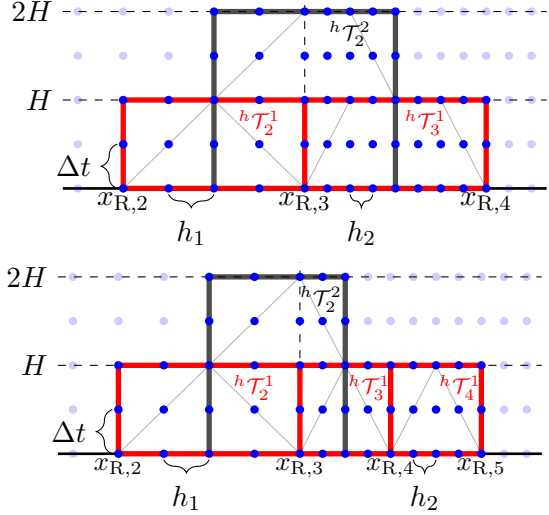
\begin{figure}
    \centering \begin{tikzpicture}[xscale=0.4,yscale=0.39]

\def\H{3}
\def\TwoH{6}


\draw[dashed] (-1,\H) -- (15,\H);
\draw[dashed] (-1,\TwoH) -- (15,\TwoH);
\draw[dashed] (7,0) -- (7,6.5);

\node[left] at (-1,\H) {$H$};
\node[left] at (-1,\TwoH) {$2H$};

\draw[line width=1pt] (-1,0) -- (15,0);


\draw[line width=2pt,black,opacity=0.7] 
(4,0) -- (4,\TwoH) -- (10,\TwoH) -- (10,0);


\draw[line width=2pt,red] 
(1,0) -- (1,\H) -- (13,\H) -- (13,0);

\draw[line width=2pt,red] 
(1,0) -- (13,0);
\draw[line width=2pt,red] (7,0) -- (7,\H);


\foreach \x in {1, 2.5, 4, 5.5, 7, 7.75, 8.5, 9.25, 10, 10.75, 11.5, 12.25, 13}
{
    \foreach \y in {0,1.5,3}
    {
        \fill[blue] (\x,\y) circle (4pt);
    }
}

\foreach \x in {4, 5.5, 7, 7.75, 8.5, 9.25, 10}
{
    \foreach \y in {4.5,6}
    {
        \fill[blue] (\x,\y) circle (4pt);
    }
}

\foreach \x in {-0.5, 1, 2.5, 10.75, 11.5, 12.25, 13, 13.75, 14.5}
{
    \foreach \y in {0,1.5, 3, 4.5,6}
    {
        \fill[blue, opacity=0.2] (\x,\y) circle (4pt);
    }
}


\draw[decorate,decoration={brace,amplitude=6pt}]
(0.8,0) -- (0.8,1.5)
node[midway,left=6pt] {$\Delta t$};


\draw[decorate,decoration={brace,amplitude=5pt}]
(4,-0.2) -- (2.5,-0.2)  
node[midway,below=6pt] {$h_1$};

\draw[decorate,decoration={brace,amplitude=5pt}]
(9.25,-0.2) -- (8.5,-0.2)
node[midway,below=6pt] {$h_2$};


\node at (1,-0.5) { $x_{\mathrm{R},2}$};
\node at (7,-0.5) {$x_{\mathrm{R},3}$};
\node at (13,-0.5) {$x_{\mathrm{R},4}$};


\node at (8.5,5.25) {\scriptsize$^h\mathcal{T}_2^2$};

\node[red] at (5.5,2.25) {\scriptsize$^h\mathcal{T}_2^1$};
\node[red] at (11,2.25) {\scriptsize$^h\mathcal{T}_3^1$};

\draw[ opacity=0.3] (1,0) -- (4,3);
\draw[ opacity=0.3] (4,3) -- (7,0);

\draw[ opacity=0.3] (7,0) -- (8.5,3);

\draw[ opacity=0.3] (11.5,3) -- (13,0);

\draw[ opacity=0.3] (4,3) -- (7,6);
\draw[ opacity=0.3] (8.5,6) -- (10,3);

\end{tikzpicture}   
\begin{tikzpicture}[xscale=0.4,yscale=0.39]

\def\H{3}
\def\TwoH{6}


\draw[dashed] (-1,\H) -- (15,\H);
\draw[dashed] (-1,\TwoH) -- (15,\TwoH);
\draw[dashed] (7,0) -- (7,6.5);

\node[left] at (-1,\H) {$H$};
\node[left] at (-1,\TwoH) {$2H$};

\draw[line width=1pt] (-1,0) -- (15,0);


\draw[line width=2pt,black,opacity=0.7] 
(4,0) -- (4,\TwoH) -- (8.5,\TwoH) -- (8.5,0);


\draw[line width=2pt,red] 
(1,0) -- (1,\H) -- (13,\H) -- (13,0);

\draw[line width=2pt,red] 
(1,0) -- (13,0);
\draw[line width=2pt,red] (7,0) -- (7,\H);
\draw[line width=2pt,red] (10,0) -- (10,\H);


\foreach \x in {1, 2.5, 4, 5.5, 7, 7.75, 8.5, 9.25, 10, 10.75, 11.5, 12.25, 13}
{
    \foreach \y in {0,1.5,3}
    {
        \fill[blue] (\x,\y) circle (4pt);
    }
}

\foreach \x in {4, 5.5, 7, 7.75, 8.5}
{
    \foreach \y in {4.5,6}
    {
        \fill[blue] (\x,\y) circle (4pt);
    }
}

\foreach \x in {-0.5, 1, 2.5, 9.25, 10, 10.75, 11.5, 12.25, 13, 13.75, 14.5}
{
    \foreach \y in {0,1.5, 3, 4.5,6}
    {
        \fill[blue, opacity=0.2] (\x,\y) circle (4pt);
    }
}


\draw[decorate,decoration={brace,amplitude=6pt}]
(0.8,0) -- (0.8,1.5)
node[midway,left=6pt] {$\Delta t$};


\draw[decorate,decoration={brace,amplitude=5pt}]
(4,-0.2) -- (2.5,-0.2)  
node[midway,below=6pt] {$h_1$};

\draw[decorate,decoration={brace,amplitude=5pt}]
(11.5,-0.2) -- (10.75,-0.2)
node[midway,below=6pt] {$h_2$};


\node at (1,-0.5) {\small $x_{\mathrm{R},2}$};
\node at (7,-0.5) {\small$x_{\mathrm{R},3}$};
\node at (10,-0.5) {\small$x_{\mathrm{R},4}$};
\node at (13,-0.5) {\small$x_{\mathrm{R},5}$};


\node at (7.7,5.25) {\scriptsize $^h\mathcal{T}_2^2$};

\node[red] at (5.5,2.25) {\scriptsize$^h\mathcal{T}_2^1$};
\node[red] at (9.3,2.25) {\scriptsize $^h\mathcal{T}_3^1$};
\node[red] at (12.3,2.25) {\scriptsize$^h\mathcal{T}_4^1$};

\draw[ opacity=0.3] (1,0) -- (4,3);
\draw[ opacity=0.3] (4,3) -- (7,0);

\draw[ opacity=0.3] (7,0) -- (8.5,3);
\draw[ opacity=0.3] (8.5,3) -- (10,0);

\draw[ opacity=0.3] (10,0) -- (11.5,3);
\draw[ opacity=0.3] (11.5,3) -- (13,0);

\draw[ opacity=0.3] (4,3) -- (7,6);
\draw[ opacity=0.3] (7,6) -- (8.5,3);

\end{tikzpicture}   \caption{Discrete subdomains for $L_1=L_2$ (top) and $m_2/m_1 = L_1/L_2=c_1/c_2$ (bottom).}
    \label{fig:discrete_subs}
\end{figure}
The wave equation is discretized by the leapfrog scheme
\begin{equation*}
\medmuskip=0mu
\thinmuskip=0mu
\thickmuskip=0mu
\nulldelimiterspace=1.2pt
\scriptspace=1.2pt    
\arraycolsep1.2em
u_{j,k}^{n+1,\ell} = 2(1-\alpha_i) u_{j,k}^{n,\ell} - u_{j,k}^{n-1,\ell} + \alpha_i ( u_{j,k}^{n,\ell+1} + u_{j,k}^{n,\ell-1} )
\end{equation*}
where $n$ and $\ell$ are time and space grid indices, $\alpha_i \coloneq \frac{c_i^2 \Delta t^2}{h_i^2}, i=1,2$, is the CFL number, and $j$ and $k$ are the space-time subdomain indices.

\FloatBarrier

\begin{figure}
\includegraphics[scale=0.6]{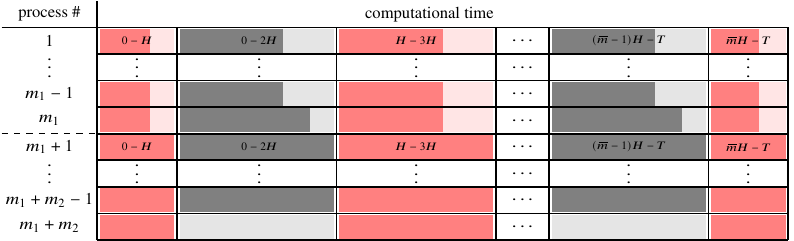} \\[2mm]
\includegraphics[scale=0.6]{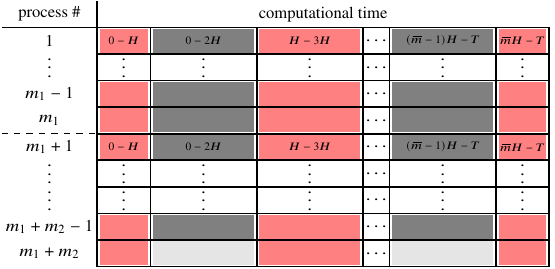}\\
 \caption{Pipelines for \(L_1=L_2\) (top) and \(m_2/m_1 = L_1/L_2 = c_1/c_2\) (bottom).}
 \label{fig:pipelines}
\end{figure}

\section{Computational cost analysis} 
In \cite{UTPhet}, we considered $m_1+m_2$ parallel processes, one per subdomain solve, and modeled the cost of solve $(j,k)$ by the subdomain measure $|\mathcal{T}_j^k|$. In that continuous framework, the decomposition choice $H_1=H_2$ and $m_1=m_2$ (equivalently, $L_1=L_2$) is optimal in terms of computational time and load balancing. At the discrete level, however, the situation changes due to the mismatch between the continuous characteristics (used in \cite{UTPhet}) and the numerical ones. The latter stem from the explicit scheme, depend on the CFL condition, and are obtained by tracking stencil-based information propagation. 

In our discrete setting, the cost of solve $(j,k)$ is modeled by $|^h\mathcal{T}_j^k|$: the number of space–time grid points in $^h\mathcal{T}_j^k$, namely the number of variables computed by the solve in the discrete subdomain $(j,k)$, see Fig.~\ref{fig:discrete_subs}.
In~\cite{UTPhet}, we showed that the height of a subdomain must not exceed the height of the solution front obtained after the first solve on that subdomain. Otherwise, additional iterations of the corresponding subproblem are required, slowing down the overall solution process. 
This remains true in the discrete setting, and we assume $H_1=H_2$ (same number of time steps within the discrete subdomains). Moreover, to avoid dissipation and dispersion in the numerical solution, we assume  $\alpha_i=1, i=1,2$, guaranteeing that numerical and continuous characteristics coincide. 
This implies that $h_1 \neq h_2$, since $c_1 >c_2$. Now, only two cases are possible: $L_1 = L_2$ and $L_1/L_2 = c_1/c_2$.
In the first case $|\mathcal{T}_j^k|=|\mathcal{T}_{i}^k|$ while $|^h\mathcal{T}_j^k|<|^h\mathcal{T}_{i}^k|$ for $j< m_1$ and $i>m_1$.
In contrast, in the second case $|\mathcal{T}_j^k| > |\mathcal{T}_{i}^k|$ while $|^h\mathcal{T}_j^k|=|^h\mathcal{T}_{i}^k|$ for $j< m_1$ and $i>m_1$.
See also Fig.~\ref{fig:discrete_subs}.
This means that, while in the continuous setting the subdomain costs are balanced for $L_1 = L_2$, the situation is opposite in the discrete setting, where the balance occurs for $L_1/L_2 = c_1/c_2$. This fact leads to the pipelines in Fig.~\ref{fig:pipelines}, where the case $L_1/L_2 = c_1/c_2$ is the optimal one (at the expense of using more processes).


\end{wavespaper}

\end{document}